\documentclass[preprint,12pt]{article}
\usepackage{natbib}
\usepackage{rotating}
\usepackage{epsfig}
\usepackage{ifpdf}
\usepackage{dsfont}
\usepackage{amssymb}

\usepackage{amssymb}

\newcommand{\X}{\mathcal{X}}

\newcommand{\R}{\mathbb{R}}

\newcommand{\E}{\mathbb{E}}
\newcommand{\C}{\mathcal{C}}

\newtheorem{thm}{Theorem}[section]
\newtheorem{cor}{Corollary}[section]
\newtheorem{pro}{Proposition}[section]
\newtheorem{rem}{Remark}[section]

\title{{\bf Large Deviation Results for the Nonparametric Regression Function Estimator on Functional Data}}

\author{\bf D. LOUANI$^a$\thanks{Corresponding author email : djamal.louani@upmc.fr}\ \ \&
S. M. OULD MAOULOUD$^b$\\ 
{\small $^a$ Universit\'e de Reims and L.S.T.A.,
Universit\'e de Paris 6, France.}\\
{\small $^b$ \'Ecole des Mines de Mauritanie, Mauritania.}}
\date{}
\begin{document}
\maketitle
\begin{quote}
\noindent {\bf Abstract} - This paper is devoted to the study of large deviation behaviors in the
setting of the estimation of the regression function on functional
data. A large deviation principle is stated for a process $Z_n$,
defined below,  allowing to derive a pointwise large deviation
principle for the Nadaraya-Watson-type $l$-indexed regression
function estimator as a by-product. Moreover, a uniform over
VC-classes Chernoff type large deviation result is stated for the
deviation of the $l$-indexed regression estimator. \vspace{1mm}

\vspace{1mm}

\end{quote}
\setcounter{section}{1}
\centerline{\large 1. INTRODUCTION}
\vspace{2mm}

\noindent The regression problem has received a great interest and
has motivated a great number of investigations and studies
throughout the time bringing to the statistic literature a
considerable knowledge.
A number of models and nonparametric estimators to estimate the
regression function have been proposed in the literature when the
discrete time or the continuous time explanatory random variables
take their values in a finite dimensional space where the Lebesgue
measure plays an important role. We refer to Bosq (1998) for an
account of properties and results and the references therein.
Due to the availability of computing resources that allow sharp
recordings in phenomena observation up to the level where data may
be treated as curves, functional modeling has received a lot of
attention in the last few years from mathematical, probabilistic,
statistical or physicist points of view. It is worth noticing that
there is an increasing number of sources of potential applications
of functional models, as in chemiometrics, environmetrics,  speech
recognition, radar range profile studies, medical data and so on.
The number of publications studying properties of these models, as
asymptotic issues for example, grows continuously. For an overview
of the present state of the art, we refer to the works of Gasser
{\it et al} (1998), Bosq (2000), Ferraty and Vieu (2000, 2004),
Ramsay and Silverman (2002, 2005),  Masry (2005), Ferraty {\it et
al} (2007), Ezzahrioui {\it et al} (2008), and to the recent
monograph by Ferraty and Vieu (2006) and the references
therein.\vspace{3mm}

\noindent To introduce the study framework, let $(X_i, Y_i)_{i\in
\mathbb{N}}$ be a sequence of i.i.d. pairs of random elements
where $Y_i$ is a real-valued random variable with density $g$,
with respect to the Lebesgue measure on $\mathbb{R}$, and $X_i$
takes its values in some semi-metric abstract space $\left( {\cal
E}, d(\cdot,\cdot)\right)$. This covers the case of semi-normed
spaces of possibly infinite dimension (e.g., Hilbert or Banach
spaces) with the norm $\Vert \cdot\Vert$ and the distance $d(x,
y)=\Vert x-y\Vert$. For a real function $l$ and any fixed $x\in
{\cal E}$, the $l$-indexed regression function at $X_1=x$ is
defined by $r^l(x) :=\mathbb{E}(l(Y_1)| X_1=x)$.\vspace{3mm}

\noindent The Nadaraya-Watson type estimator of $r^l$ has been
introduced by Ferraty and Vieu (2000). It is defined, for any
fixed $x\in{\cal E}$, by
\begin{eqnarray}\label{EstimateRgression.1}
\hat{r}_{n}^l(x) &=& \left\{ \begin{array}{ll} \frac{ \sum_{i=1}^n
l(Y_i) K\left( \frac{d(x, X_i)}{h} \right)} { \sum_{i=1}^n K\left(
\frac{d(x, X_i)}{h}
\right)}:=\frac{\hat{r}_{n,2}(x)}{\hat{r}_{n,1}(x)}, & if\ \
\hat{r}_{n,1}(x)\neq 0,\\
0,& elsewhere.
\end{array}
\right.
\end{eqnarray}
Here, $K$ is a real-valued kernel function, $h:=h_n$ is the
bandwidth parameter (which goes to $0$ as $n$ goes to infinity),
\begin{eqnarray}\label{EstimateRgression.1.1}
\hat{r}_{n,1}(x)=\frac{1}{n\phi(h)}\sum_{i=1}^n \Delta_i(x) \ \
\mbox{and}\ \ \hat{r}_{n,2}(x)=\frac{1}{n\phi(h)}\sum_{i=1}^n
l(Y_i)\Delta_i(x),
\end{eqnarray}
where
$$\Delta_i(x)=K\left( \frac{d(x, X_i)}{h}\right)$$
and $\phi$ is a positive function that will be defined below.
Notice that the index function $l$ allows to study simultaneously
properties of several estimates. The first example is given by the
most classical regression function estimator where $l$ stands as
the identity function. Whenever $l=\mathds{1}_A$ is the indicator
function of the set A, $\hat{r}_{n}^l(x)$ is the estimator of the
conditional probability measure of the event $\{Y_1\in A\}$ given
$X_1=x$.\vspace{3mm}

\noindent For $x\in {\cal E}$, consider now the vector process
$$Z_n(x)=(\hat{r}_{n,1}(x),\hat{r}_{n,2}(x)).$$
In this paper, we aim at establishing a large deviation principle for
the process $Z_n(x)$ and deriving asymptotics for the $l$-indexed
regression function estimator in both the pointwise case and the uniform, over some vc-classes, case.\vspace{3mm}

\noindent There exists an extensive large deviation literature
involving many areas of probability and statistics. We refer to
the books of Dembo \& Zeitouni (1998) and Deuschel \& Stroock
(1989) and the references therein for an account of results and
applications. In nonparametric function estimation setting,
several results have been obtained these last years. We refer to
Louani (1999) and Ould Maouloud (2008) where the studies involve
the Nadaraya-Watson and histogram estimates of the regression
function respectively both in the real vector case. Using the
delta-sequence estimation method, Louani and Ould Maouloud (2011)
established a large deviation principle for the real regression
function estimate embedded in the $L_1$ space equipped with the
weak topology. Notice that the main applications of large
deviation results are related to the efficiency of tests in the
Bahadur sense, see Nikitin (1995) for more details,  together with
the inaccuracy rate of estimators that allow to compare testing
procedures and estimation performances respectively. The results
may be also used to establish estimates consistency with rates of
convergence.\vspace{5mm}

\setcounter{section}{2}
\centerline{\large 2. RESULTS}
\vspace{2mm}

\noindent Our results are stated under some  assumptions we
gather hereafter for easy reference
\begin{enumerate}
\item[(A1)] $K$ is a nonnegative bounded differentiable kernel
over its support $[0, 1]$ and $K(1)>0$. The derivative $K^\prime$
of $K$  exists on the interval $[0, 1]$.

\item[(A2)] For $x\in{\cal E}$ and a real number $v$, there exist
a nonnegative functional $f_{v}$ and a nonnegative real function
$\phi$ tending to zero, as its argument tends to $0$, such that,
uniformly in $v$,

(i) $F_x(u|Y=v)=\mathbb{P}(d(x,X_1)\leq u|Y=
v)=\phi(u)f_v(x)+o(\phi(u))$ as $u\to 0$,

(ii) There exists a nondecreasing bounded function $\tau_0$ such
that, uniformly in $u\in [0, 1]$,
$$\frac{\phi(hu )}{\phi(h)}
=\tau_0(u)+o(1),\ \ \mbox{as}\ \ h\downarrow 0,$$

\item[(A3)] 
For any real numbers $a$ and $b$,
\begin{eqnarray*}
\hspace*{-2cm}&&(i)\ \int(e^{a+bl(v))}-1)f_v(x)g(v)dv<\infty,\ \ (ii)\ \int e^{al(v)}f_v(x)g(v)dv<\infty,  \\
\hspace*{-2cm}&&(iii)\ \int e^{al(v)}f_v(x)g(v)l(v)dv<\infty,\ \
(iv)\ \int e^{al(v)}f_v(x)g(v)l^2(v)dv<\infty.
\end{eqnarray*}
\end{enumerate}

\noindent{\bf Discussions of hypotheses.}  Condition (A1) is  very
usual in nonparametric estimation literature devoted to functional
data context. From the fact that Lebesgue measure does not exist
on infinite dimension space, hypotheses (A2) involve the small
ball techniques related to the fractal dimension used in this
paper. A number of examples of the function $\phi$ together with
the corresponding decomposition of the probability of the small
balls are given throughout several works (See, e.g., Ferraty and
Vieu (2000, 2004, 2006), Ferraty {\it et al} (2006), Ezzahrioui
{\it et al} (2008) and La\"ib and Louani (2010)). A further
example is given hereafter to illustrate the condition (A2)(i).
Hypotheses (A3) are set on to insure the needed properties of
finiteness and differentiability of the moment generating function
of the process $Z_n(x)$. These hypotheses induce the fact that the
large deviation principle holds with a good rate function, a
property
that is strongly expected in such results. 
\vspace{5mm}

\noindent{\bf Example 1.}\ Let ${\cal E}$ be a separate abstract
space equipped with the semi-metric defined, for $(x, y)\in {\cal
E}^2$, by
$$d(x, y)=|\int (x(t)-y(t))dt|.$$
Consider two elements $h$ and $l$ in ${\cal C}$ together with the
regression model
$$X_i=Y_{i}h+\varepsilon_il,$$
where $\varepsilon_i$ is a real random variable independent of
$Y_i$.
Observe now, for any $u>0$, that we have
$$F_x(u|Y=v)=\mathbb{P}\left(d(x, X_i) \leq u | Y_i=v\right)=
\mathbb{P}\left(|\int(x(t)-X_i(t))dt| \leq u | Y_i=v\right).$$
Consequently, while assuming $0\neq |\int l(t)dt|<\infty$, $|\int
x(t)dt|<\infty$ and
 $|\int h(t)dt|<\infty$, it follows that
\begin{eqnarray*}
F_{x}(u|Y=v)&=& \mathbb{P}\left(|\int(x(t)-Y_ih(t)-\varepsilon_il(t))dt| \leq u | Y_{i}=v\right)\\
&=& \mathbb{P}\left(|\int(x(t)-vh(t)-\varepsilon_il(t))dt| \leq u \right)\\
&=&
\mathbb{P}\left(\frac{-u+\int x(t)dt-v\int h(t)dt}{\int l(t)dt} \leq \varepsilon_i\leq \frac{u+\int x(t)dt-v\int h(t)dt}{\int l(t)dt}\right) \\
&=& \Phi\left( \frac{u+\int x(t)dt-v\int h(t)dt}{\int l(t)dt}
\right)- \Phi\left(\frac{-u+\int x(t)dt-v\int h(t)dt}{\int
l(t)dt}\right),
\end{eqnarray*}
where $\Phi$ is the distribution function of $\varepsilon_i$.
Taking $\Phi$ as the ${\cal N}(0, 1)$ distribution function and
assuming that $0<\int l(t)dt<\infty$, we obtain,
\begin{eqnarray*}
F_x(u|Y=v) &= & \frac{u}{\int
l(t)dt}\sqrt{\frac{2}{\pi}}\exp\left( -\frac{1}{2}\left(
\frac{\int x(t)dt- v\int h(t)dt}{\int l(t)dt} \right)^2\right)
(1+o(1)),
\end{eqnarray*}
and the condition (A2)(i) is satisfied with $\phi(u)=2u$ and
$$f_v(x)=\frac{1}{\int l(t)dt\sqrt{2\pi}}\exp\left(
-\frac{1}{2}\left( \frac{\int x(t)dt- v\int h(t)dt}{\int l(t)dt}
\right)^2\right).$$
\ \ \ \ \ \ \ \ \ \ \  \ \ \ $\hfill\Box$

\vspace{5mm}

\noindent The first result states a LDP for the process $Z_n(x)$.

\begin{thm}\label{thm1}
Under assumptions (A1)-(A3), $Z_n(x)$ satisfies a LDP with the
speed $n\phi(h)$ and a good rate function given by
$$\Gamma_x(\lambda_1,\lambda_2)=\sup_{t_1,t_2}\left\{\lambda_1t_1+\lambda_2t_2-\Phi^x(t_1,t_2)\right\},$$
where
\begin{eqnarray*}
\Phi^x(t_1,t_2)&=&\int f_v(x)\left[\left(e^{(t_1+t_2l(v))K(1)}-1\right)\right.\\
&&-\left.\int_0^1(t_1+t_2l(v))K'(u)e^{(t_1+t_2l(v))K(u)}\tau_0(u)du\right]g(
v)dv.
\end{eqnarray*}
\end{thm}
\begin{rem} If we suppose that the function $\tau_0$ is
differentiable, then, integrating by parts,  we obtain
$$\Phi^x(t_1,t_2)=\int\int_0^1f_v(x)\tau^\prime_0(u)\left(e^{(t_1+t_2l(v))K(u)}-1\right)g(v)dudv$$
which gives a more simpler form of the rate function.
\end{rem}
\vspace{3mm}

\noindent Whenever we take the function $K$ as the uniform kernel,
we obtain a more explicit rate function. In order to display it,
we introduce first some further notations. For any
$x\in\mathcal{E}$ and any $t\in\mathbb{R}$, set
\begin{equation}\label{defV}
V_x(t)=\frac{\int e^{tl(u)}f_u(x)g(u)l(u)du}{\int e^{tl(v)}f_v(x)g(v)dv}
\end{equation}
and
$$
V_x^{-1}(t)=\inf\{s : V_x(s)\geq t\}.
$$
Moreover, assuming that the derivative of the function $\tau_0$
exists and considering the fact that $\tau_0(0)=0$ and
$\tau_0(1)=1$, we observe, whenever $K(u)=\mathds{1}_{[0,1]}(u)$,
that
\begin{equation}\label{kerunif}
\Phi^x(t_1,t_2)=\int \left(e^{(t_1+t_2l(v))}-1\right)f_v(x)g(v)dv.
\end{equation}

\begin{cor}\label{cor1}
Assume that the function $\tau_0$ is differentiable and that $K$
is the uniform kernel over the interval $[0,1]$, then, under
assumptions (A3)(i)-(iv), we obtain the following explicit form of
the rate function
\begin{eqnarray}\label{rfunct}
\Gamma_x(\lambda_1,\lambda_2)\!\!=\!\!\left\{\begin{array}{l}
\lambda_1(\log\lambda_1
-1)+\lambda_2V_x^{-1}\left(\frac{\lambda_2}{\lambda_1}\right)-
\lambda_1\log\int
e^{V_x^{-1}\left(\frac{\lambda_2}{\lambda_1}\right)l(v)}f_v(x)g(v)dv\\
+\int f_u(x)g(u)du \ \ \  if\ \ \lambda_1>0\ \ \ and\ \ \ \ v_0(x)<\lambda_2/\lambda_1 <v_1(x)\\
\infty \ \ \ \ \ \ \ \ \ \ \ \  \ \ \  \ \ \  \ \  elsewhere,
\end{array}
\right.
\end{eqnarray}
where $v_0(x)=\inf_{t\in\mathbb{R}}V_x(t)$ and
$v_1(x)=\sup_{t\in\mathbb{R}}V_x(t)$.
\end{cor}
\vspace{3mm}

\begin{rem}
 Whenever $l=\mathds{1}_A$ is the indicator function of a subset $A$ of $\mathbb{R}$, it is possible to display a more explicit rate function
 whenever $\tau_0$ is differentiable. Towards this
 end, for any $B\subset\mathbb{R}$ and any $t\in\mathbb{R}$, set
$$W_x(B)=\int_Bf_v(x)g(v)dv\ \ \mbox{and}\ \
\zeta(t)=\int_0^1\tau_0^\prime(u)K(u)e^{tK(u)}du.$$ It follows
then that
\begin{eqnarray*}
\Gamma_x(\lambda_1,\lambda_2)&=&(\lambda_1-\lambda_2)\zeta^{-1}\left(\frac{\lambda_1-\lambda_2}{W_x(\bar{A})}\right)+\lambda_2
\zeta^{-1}\left(\frac{\lambda_2}{W_x(A)}\right)+\int f_v(x)g(v)dv\\
&&-\int_0^1\tau_0^\prime(u)\left[W_x(A)\exp\left\{\zeta^{-1}\left(\frac{\lambda_2}{W_x(A)}\right)K(u)\right\}\right.\\
&&\hspace{2cm}+\left.W_x(\bar{A})\exp\left\{\zeta^{-1}\left(\frac{\lambda_1-\lambda_2}{W_x(\bar{A})}\right)K(u)\right\}\right]du,
\end{eqnarray*}
where $\bar{A}$ is the complementary set of $A$ and
$$\zeta^{-1}(t)=\inf\{s : \zeta(s)\geq t\}.$$
\end{rem}

\noindent The following corollary gives the result pertaining to a
large deviation principe for the regression function estimate at
the point $x$.

\begin{cor}\label{cor2}
Under hypotheses of Theorem \ref{thm1}, the regression function
estimate $\hat{r}_n^l(x)$ satisfies a LDP with the speed
$n\phi(h)$ and the good rate function defined by
$$\gamma_x(\lambda)=\inf_{\lambda_1}\{\Gamma_x(\lambda_1,\lambda\times\lambda_1)\}.$$
\end{cor}

\begin{rem}\label{rem3}
If we assume that the function $\tau_0$ is differentiable and that
$K$ is the uniform kernel over the interval $[0,1]$, then we
obtain the following explicit form of the rate function
\begin{eqnarray}\label{gamma_unif}
 \gamma_x(\lambda) &=& \int\left(
1-\exp\left\{V_x^{-1}\left(\lambda\right)(l(v)-\lambda\right\}\right)
f_v(x)g(v)dv,
\end{eqnarray}
whenever $v_0(x)<\lambda<v_1(x)$, and $\gamma_x(\lambda)=\infty$
elsewhere.
\end{rem}
\vspace{1mm}
\begin{rem}\label{remark4}
The first and second derivatives of the function $\gamma_x$ given
in the statement (\ref{gamma_unif}) are
$$\gamma_x^\prime(\lambda)=V_x^{-1}(\lambda)\exp\{-\lambda
V_x^{-1}(\lambda)\}\int e^{V_x^{-1}(\lambda)l(v)}f_v(x)g(v)dv$$
and
$$\gamma_x^{\prime\prime}(\lambda)=\left(\frac{1}{V_x^\prime(V_x^{-1}(\lambda))}-(V_x^{-1}(\lambda))^2\right)\exp\{-\lambda
V_x^{-1}(\lambda)\}\int e^{V_x^{-1}(\lambda)l(v)}f_v(x)g(v)dv$$
respectively. When $Z$ denotes the random variable associated to
the density function $f_v(x)g(v)/\int f_v(x)g(v)dv$, it follows
that $V_x(0)= \E(l(Z))$. Therefore, assuming that $\E(l(Z))=0$, by Taylor series expansion we
obtain, in the neighborhood of $\lambda=0$, that
$$\gamma_x(\lambda)=\frac{\lambda^2}{2\ \E(l^2(Z))}\int f_v(x)g(v)dv (1+o(1)).$$
\end{rem}
\vspace{3mm}

\noindent In the sequel, we investigate the uniform aspects of
large deviation, in the Chernoff sense, of the regression function
estimate $\hat{r}_n^l(x)$. More precisely, we consider the
asymptotic behavior of the quantity
$\displaystyle{\|\hat{r}_n^l-r^l\|_{\cal C}:=\sup_{x\in{\cal
C}}|\hat{r}_n^l(x)-r^l(x)|}$, where ${\cal C}$ is a class of
elements of ${\cal E}$. Towards this end, for any $\varepsilon>0,$
consider the following number
\begin{eqnarray*}
{\cal N}(\varepsilon,{\cal C},d)&=&\min\{n:\ \mbox{there exist}\ c_1\cdots,c_n \ \mbox{in}\ {\cal C}\ \mbox{such that}\ \forall\ x\in {\cal C}  \\
&& \hspace{1.3cm}\ \mbox{ there exists}\ 1\leq k\leq n \
\mbox{such that}\ d(x,c_k)<\varepsilon\}
\end{eqnarray*}
which measures how full is the class ${\cal C}$. Further notations
are needed to display the uniform large deviation result. From now
on, set
\begin{eqnarray}\label{beta}
\beta(x,\lambda)=\inf\{\gamma_x(\alpha+r^l(x)) : \alpha\in
(-\infty,-\lambda]\cup[\lambda,\infty)\},
\end{eqnarray}
$\rho(\lambda)=\inf_{x\in{\cal C}}\beta(x,\lambda)$,\
 and $D_Y$ to be the set of values
of the random variable $Y_1$. Moreover consider the conditions
\begin{enumerate}
\item[(A4)] (i) $\displaystyle{\sup_{x\in{\cal C}}\int
|l(v)|f_v(x)g(v)dv <\infty}$,

(ii) The condition (A2)(i) is satisfied uniformly in
$\displaystyle{v\in D_Y}$,

(iii) $\displaystyle{\int_0^1|K^\prime(u)|\tau_0(u)du<\infty}$.
\end{enumerate}
\vspace{1mm}

\begin{rem}
In the setting of Remark \ref{rem3}, it is easily seen that
$$\beta(x,\lambda)=\min\{\gamma_x(-\lambda+r^l(x)),\gamma_x(\lambda+r^l(x))\},$$
since the function $\gamma_x$ is non-increasing on the left of
$r^l(x)$ and non-decreasing on its right.
\end{rem}
\vspace{1mm}

\noindent The following theorem gives a Chernoff-type large
deviation result for the uniform deviation of the estimate
$\hat{r}_n^l$ with respect to $r^l$.

\begin{thm}\label{thm2}
Suppose that the function $r^l$ is uniformly continuous upon
${\cal C}$ and that the kernel $K$ is a Lipschitz function bounded
from below by a constant $K_0>0$. Under hypotheses (A1)-(A4)
whenever  the condition
\begin{equation}\label{entropy}
 \nu=o\left(\frac{nh}{\exp\{An\phi(h)\}}\right), \ \ \mbox{for any}\ \ A>0,\ \ \mbox{and}\ \ \lim_{n\rightarrow\infty}\frac{\log {\cal N}(\nu,{\cal
   C},d)}{n\phi(h)}=0
\end{equation}
is satisfied and the function $\rho$ is continuous, for any
$\lambda>0$, we have
\begin{equation}
\lim_{n\rightarrow\infty}\frac{1}{n\phi(h)}\log\mathbb{P}(\|\hat{r}_n^l-r^l\|_{\cal
C}>\lambda)=-\rho(\lambda).
\end{equation}
\end{thm}
\vspace{3mm}

\noindent The continuity of the rate function $\rho$ is a needed
condition to obtain the result of Theorem \ref{thm2}. It is then
natural to ask the question about the required assumptions for
this condition to be satisfied. The following propositions give a
reply to this question.

\begin{pro}\label{prop1}
Assume that the parametric family of functions
$\{\Gamma_x(\lambda,\mu)\}_{x\in \C}$ is  equi-continuous and that
the function $r^l(x)$ is bounded. Then, the function $\rho$ is
continuous.
\end{pro}
\vspace{3mm}

\noindent It is difficult to state the conditions under which the
family of functions $\{\Gamma_x(\lambda,\mu)\}_{x\in \C}$ is
equi-continuous in the general framework. Hereafter, we limit
ourselves to the case where the kernel $K$ is uniform over the
interval $[0,1]$ and the function $\tau_0$ is differentiable as in
the setting of  Corollary \ref{cor1}. Towards this end, we first
introduce the following  notations. From now on, $\bar{S}$ stands
as the generic notation of the complementary of any set $S$,
$\displaystyle{d(v)=\inf_{x\in\C}f_v(x)}$,
$\displaystyle{D(v)=\sup_{x\in\C}f_v(x)}$, $A_l=\{v : l(v)-1\leq
0\}$, $B_l=\{v : l(v)+1\leq 0\}$ and $C_l=\{v : l(v)\leq 0\}$.
Furthermore, we  consider these additional assumptions
\begin{enumerate}
\item[(A5)] (i) Whenever $S$ stands as one of the sets $A_l$ or
$B_l$, we have $0<\int_Sd(v)g(v)dv<\infty$ and
$0<\int_{\bar{S}}d(v)g(v)dv<\infty$,

(ii) For any $t$, $\displaystyle{\left|\int
e^{tl(v)}l(v)(D(v)\mathds{1}_{C_l}(v)+d(v)\mathds{1}_{\bar{C}_l}(v))g(v)dv\right|<\infty},$

(iii) For any $t$, $\displaystyle{\left|\int
e^{tl(v)}l(v)(d(v)\mathds{1}_{C_l}(v)+D(v)\mathds{1}_{\bar{C}_l}(v))g(v)dv\right|<\infty},$

(iv) $\displaystyle{\int_{C_l} e^{tl(v)}\left(\mathds{1}_{\{t\leq
0\}}\mathds{1}_{C_l}(v)+\mathds{1}_{\{t>
0\}}\mathds{1}_{\bar{C}_l}(v)\right)D(v)g(v)dv<\infty}$.
\end{enumerate}
\vspace{2mm}

\begin{pro}\label{prop2}
In the setting of Corollary \ref{cor1}, assume that assumptions
(A5) are satisfied. Then the family of functions
$\{\Gamma_x(\lambda,\mu)\}_{x\in \C}$ is equi-continuous on its
finiteness domain.
\end{pro}
\vspace{2mm}

\begin{rem}
Whenever the smoothing parameter $h$ and the function $\phi$ are
such that
\begin{equation}\label{phih}
\frac{\exp\{An\phi(h)\}}{n\phi(h)}=o(nh),\ \ \mbox{for any}\ \
A>0,\
\end{equation}
then the condition (\ref{entropy}) takes tha form
\begin{equation}\label{entopy-cond}
\lim_{\nu\rightarrow 0}\nu\log {\cal N}(\nu,{\cal C},d)=0.
\end{equation}
Notice that the condition (\ref{phih}) is satisfied when, for
example, we take $\phi(h)=ah^\alpha$ and $h=\left(\frac{\log\log
n}{n}\right)^{\frac{1}{\alpha}}$ with $\alpha>1$ and $a>0$. The
condition (\ref{entopy-cond}) is very usual in defining
Vapnik-Chervonenkis classes, see, for instance van der Vaart and
Wellner (1996). Hereafter, examples of classes fulfilling the
condition (\ref{entopy-cond}) are displayed.
\end{rem}
\vspace {3mm}

\noindent{\sc Example 1.}\ ({\it Parametric classes of functions})
\ \ For a function $\X$ in the space $L_p(\R)$, consider the
parametric class of functions defined by
$$\C=\{\X_{a}(.)=a\X(a.) : -\infty<A_1\leq a\leq A_2<\infty\
\mbox{and}\ a\neq 0\}$$
together with the $L_p$-distance given, for $1\leq p<\infty$, by
$$d_p(x,y)=\left(\int|x(t)-y(t)|^pdt\right)^{\frac{1}{p}}.$$
assuming that  the function $\X$ is differentiable with a
continuous derivative and that the function $q(t)=t\X^\prime(t)$
is $L_p$-integrable, it follows that there exists a positive
constant $C_0$ such that $N(\nu,\C,d_p)\leq C_0/\nu$. Therefore, the
condition (\ref{entopy-cond}) is satisfied.
 \vspace {4mm}

\noindent{\sc Example 2.}\ ({\it Classes of functions that are
Lipschitz in a parameter})\ \ Let $T$ be an index set and consider
$d$ a distance over $T$. Suppose that ${\cal C}=\{x_t : t\in T\}$
is a class of functions defined on $\mathbb{R}$, that are
Lipschitz in the index parameter $t$ in the sense that there
exists a function $x$ on $\mathbb{R}$ such that for any $u\in
\mathbb{R}$,
$$|x_s(u)-x_t(u)|\leq d(s,t)x(u).$$
From Theorem 2.7.11 in van der Vaart \& Wellner (1996), it
follows, for any norm $\|.\|$ whenever $\|x\|<\infty$, that
$${\cal N}(\nu\|x\|,{\cal C},\|.\|)\leq N(\nu,T, d),$$
where $N(\nu,T, d)$ is the minimal number of balls of radius $\nu$
needed to cover $T$. Therefore, the condition (\ref{entopy-cond})
may be expressed as
$$\lim_{\nu\rightarrow 0}\nu\log N(\nu,T,d)=0,$$
which is a Vapnik-Chervonenkis class of sets condition. Naturally,
it is more easy to display examples of Vapnik-Chervonenkis classes
of sets. As an example, let $x$ be a Lipschitz boundedly supported
function. For an index set $T$ included in $\mathbb{R}$, define
the class ${\cal C}$ by taking, for any $u\in\mathbb{R},
x_t(u)=x(u-t)$. It is obvious then that the $x_t$'s are Lipschitz
with respect to the index parameter $t$. Taking $d$ as the
absolute distance on $T$ that we take as a bounded convex
interval, it follows easily that $N(\nu,T,d)\leq |T|/\nu$, where
$|T|$ is the diameter of $T$. Therefore, the condition
(\ref{entopy-cond}) is satisfied.

\vspace{4mm}

\noindent{\sc Example 3.}\ ({\it Smooth function classes})\ \  Let
${\cal U}$ be a bounded convex subset of $\mathbb{R}^d$ with
nonempty interior. For any $\alpha>0$, consider the class of
functions on ${\cal U}$ that possess uniformly bounded derivatives
up to order $[\alpha]$, where $[\alpha]$ stands as the integer
part of $\alpha$, and whose highest derivatives are Lipschitz of
order $\alpha-[\alpha]$. Denote by $x^{(l)}$ the $l$-th derivative
of $x$ and set, for any $x : {\cal U}\rightarrow \mathbb{R}$,
$$\|x\|_{\alpha}:=\max_{l\leq[\alpha]}\sup_{u\in{\cal U}}|x^{(l)}(u)|+
\sup_{{u,v\in{\stackrel{o}{{\cal U}}},u\neq v}}
\frac{|x^{[\alpha]}(u)-x^{[\alpha]}(v)|}{|u-v|^{\alpha-[\alpha]}}$$
where ${\stackrel{o}{{\cal U}}}$ is the interior of ${\cal U}$.
Let ${\cal C}=C_M^{\alpha}({\cal U})$ be the class of all
continuous functions $x : {\cal U}\rightarrow \mathbb{R}$ with
$\|x\|_{\alpha}\leq M$. It follows from Corollary 2.7.2 of van der
Vaart \& Wellner (1996) that, for every $p\geq 1$ and any $\nu>0$,
there exists a constant $c$ depending only on $\alpha$, the
diameter of ${\cal U}$ and $d$ such that
$$\log N_{[\ ]}(\nu,C_1^{\alpha}({\cal U}),L_p)\leq c\left(\frac{1}{\nu}\right)^{\frac{d}{\alpha}}.$$
Here, ${\cal N}_{[\ ]}$ denotes the bracketing number (see, for
instance, van der Vaart \& Wellner (1996) page 83, for the
definition) and, for a probability measure $\mu$,
$L_p(x,y)=(\int|x(u)-y(u)|^pd\mu(u))^{\frac{1}{p}}$. Since
$$N(\nu,C_1^{\alpha}({\cal U}),L_p)\leq N_{[\
]}(2\nu,C_1^{\alpha}({\cal U}),L_p),$$ it is then clear that the
condition (\ref{entopy-cond}) is satisfied provided that
$\alpha>1$.\vspace{5mm}

\setcounter{section}{2}
\centerline{\large 3. PROOFS}
\vspace{2mm}

\noindent{\bf Proof of Theorem \ref{thm1}}\ \
The Laplace transform associated to the process $n\phi(h)Z_n(x)$
is defined, for any $(t_1,t_2)\in\mathbb{R}^2$, by
\begin{eqnarray*}
\Phi_n^x(t_1,t_2)&=&\mathbb{E}\left[\exp\{<(t_1,t_2),n\phi(h)(\hat{r}_{n,1}(x),\hat{r}_{n,2}(x))>\}\right]\\
&=&\mathbb{E}\left[\exp\left\{\left<(t_1,t_2),\left(\sum_{i=1}^n\Delta_i(x),\sum_{i=1}^nl(Y_i)\Delta_i(x)\right)\right>\right\}\right]\\
&=&\mathbb{E}\left[\exp\left\{\sum_{i=1}^n(t_1+t_2l(Y_i))\Delta_i(x)\right\}\right]\\
&=&\left(\mathbb{E}\left[\exp\{(t_1+t_2l(Y_i))\Delta_1(x)\}\right]\right)^n:=\left(\varphi_n^x(t_1,t_2)\right)^n,
\end{eqnarray*}
where $<.,.>$ denotes the inner product. Let us now evaluate the
quantity $\varphi_n^x(t_1,t_2)$. Observe that
\begin{eqnarray*}
\varphi_n^x(t_1,t_2)&=&1+\mathbb{E}\left[\exp\{(t_1+t_2l(Y_1))\Delta_1(x)\}-1\right]\\
&=&1+\int_0^1\int\left(e^{(t_1+t_2l(v))K(u)}-1\right)d\mathbb{P}\left(\frac{d(x,X_1)}{h}\leq
u,Y\leq v\right)\\
&=&1+\int_0^1\int\left(e^{(t_1+t_2l(v))K(u)}-1\right)d\mathbb{P}\left(\frac{d(x,X_1)}{h}\leq
u\left|\right.Y=v\right)d\mathbb{P}(Y\leq v)
\end{eqnarray*}
Integrating by parts with respect to the component $u$, we obtain
\begin{eqnarray*}
\varphi_n^x(t_1,t_2)&=&1+\int\left[\left(e^{(t_1+t_2l(v))K(1)}-1\right)F_x(h|Y=v)\right.\\
&&\ \ \
-\int_0^1\left.(t_1+t_2l(v))K'(u)e^{(t_1+t_2l(v))K(u)}F_x(uh|Y=v)du\right]d\mathbb{P}(Y\leq
v).
\end{eqnarray*}
Making use of the condition (A2)(i), we obtain
\begin{eqnarray*}
\varphi_n^x(t_1,t_2)&=&1+\int\left[\left(e^{(t_1+t_2l(v))K(1)}-1\right)(\phi(h)f_v(x)+o(\phi(h)))\right.\\
&-&\int_0^1\left.(t_1+t_2l(v))K'(u)e^{(t_1+t_2l(v))K(u)}(\phi(uh)f_v(x)+o(\phi(uh))du\right]d\mathbb{P}(Y\leq
v)\\
&=&1+\phi(h)\left(\int\left[\left(e^{(t_1+t_2l(v))K(1)}-1\right)(f_v(x)+o(1))\right.\right.\\
&-&\int_0^1\left.\left.(t_1+t_2v)K'(u)e^{(t_1+t_2l(v))K(u)}(f_v(x)+o(1))\frac{\phi(uh)}{\phi(h)}du\right]d\mathbb{P}(Y\leq
v)\right).
\end{eqnarray*}
By the condition (A2)(ii), it follows that
\begin{eqnarray*}
\varphi_n^x(t_1,t_2)\!\!\!&=&1+\phi(h)\left(\int\left[\left(e^{(t_1+t_2l(v))K(1)}-1\right)(f_v(x)+o(1))\right.\right.\\
&-&\!\!\!\!\int_0^1\!\!\left.\left.(t_1\!+\!t_2l(v))K'(u)e^{(t_1+t_2l(v))K(u)}(f_v(x)\!+\!o(1))(\tau_0(u)\!+\!o(1))du\right]d\mathbb{P}(Y\leq
v)\right).
\end{eqnarray*}
Therefore, after a Taylor series expansion of the function
$\log(1+u)$ around $u=0$, we obtain
\begin{eqnarray*}
\lim_{n\rightarrow\infty}\frac{1}{n\phi(h)}\log\Phi_n^x(t_1,t_2)&:=&\Phi^x(t_1,t_2)=\int
f_v(x)\left[\left(e^{(t_1+t_2l(v))K(1)}-1\right)\right.\\
&-&\int_0^1\left.(t_1+t_2l(v))K'(u)e^{(t_1+t_2l(v))K(u)}\tau_0(u)du\right]g(
v)dv
\end{eqnarray*}
Note that the condition (A3) implies that the function
$\Phi^x(t_1,t_2)$ is finite and  differentiable everywhere. The
Fenchel-Legendre transform of $\Phi^x(t_1,t_2)$ is given by
\begin{eqnarray*}
\Gamma_x(\lambda_1,\lambda_2)=\sup_{t_1,t_2}\left\{\lambda_1t_1+\lambda_2t_2-\Phi^x(t_1,t_2)\right\}.
\end{eqnarray*}
We have now to establish that the function $\Phi^x(t_1,t_2)$ is
essentially smooth and to use the G\"artner-Ellis Theorem (see,
Dembo \& Zeitouni (1998), page 44) to achieve the proof.

Considering hypotheses (A3) (i)-(ii), it is clear that the
interior of the set $D=\{(t_1,t_2) : \Phi^x(t_1,t_2)<\infty\}$ is
not empty. Moreover, making use of conditions (A3), it follows
that the function $\Phi^x(t_1,t_2)$ is differentiable throughout
the domain $\stackrel{\circ}{D}$. Subsequently, it is clear  that
the function $\Phi^x(t_1,t_2)$ is steep and, therefore, is
essentially smooth. \hfill{$\Box$}\vspace{5mm}

\noindent{\bf Proof of Corollary \ref{cor1}}\ \ In view of the
statement (\ref{kerunif}), we have to maximize the function
$$Q(t_1,t_2)=\lambda_1t_1+\lambda_2t_2-\int
\left(e^{(t_1+t_2l(v))}-1\right)f_v(x)g(v)dv.$$ Since the function
$Q$ is concave, it is easily seen that its maximum is reached at
the point
$$(t_1,t_2)=\left(\log(\lambda_1)-\log\int
e^{V_x^{-1}\left(\frac{\lambda_2}{\lambda_1}\right)l(v)}f_v(x)g(v)dv,V_x^{-1}\left(\frac{\lambda_2}{\lambda_1}\right)\right)$$
which gives the main form of the rate function given in the
statement (\ref{rfunct}).

In order to display the finiteness domain of the function
$\Gamma_x$, we have to study the function $V_x$. Observe from
hypotheses (A3) that $V_x$ is a differentiable function and that
its derivative is given by
\begin{eqnarray*}
V_x^\prime(t)&=&\frac{\left(\int
e^{tl(v)}f_v(x)g(v)l^2(v)dv\right)\left(\int e^{tl(v)}
f_v(x)g(v)dv\right)- \left(\int
e^{tl(v)}f_v(x)g(v)l(v)dv\right)^2}{\left(\int
e^{tl(v)}f_v(x)g(v)dv\right)^2}\\
&=&\frac{\int\int
e^{t(l(u)+l(v))}f_u(x)f_v(x)g(u)g(v)(l^2(v)-l(u)l(v))dudv}{\left(\int
e^{tl(v)}f_v(x)g(v)dv\right)^2}\\
&=&\frac{1}{2}\frac{\int\int
e^{t(l(u)+l(v))}f_u(x)f_v(x)g(u)g(v)(l(v)-l(u))^2dudv}{\left(\int
e^{tl(v)}f_v(x)g(v)dv\right)^2}\geq 0.
\end{eqnarray*}
Therefore,  $V_x$ is an increasing function. Notice that
$v_1(x)=\lim_{t\rightarrow\infty}V_x(t)=\sup_tV_x(t)$ exists in
the closure $\bar{\mathbb{R}}$ of $\mathbb{R}$.

Assuming now that $\lambda_2/\lambda_1>v_1(x)$, it results that
there exists $\varepsilon>0$ such that, for any $t\in\mathbb{R}$,
\begin{equation}\label{domain}
\frac{\lambda_2}{\lambda_1}+\varepsilon\geq V_x(t).
\end{equation}
Integrating in both sides of (\ref{domain}) with respect to $t$,
it follows, for any $t\in\mathbb{R}$, that
$$\left(\frac{\lambda_2}{\lambda_1}+\varepsilon\right)t+c_0\geq \log\int
e^{tl(v)}f_v(x)g(v)dv,$$ with $c_0=\log\int f_v(x)g(v)dv$.
Therefore, for any $(t_1,t_2)\in\mathbb{R}^2$, we obtain
\begin{eqnarray*}
  H(t_1,t_2) &:=& \lambda_1t_1+\lambda_2t_2-\exp\left\{t_1+ \left(\frac{\lambda_2}{\lambda_1}+\varepsilon\right)t+c_0\right\}+e^{c_0}\\
  &\leq&
  \lambda_1t_1+\lambda_2t_2-\int\left(\exp\left\{t_1+t_2l(v)\right\}-1\right)f_v(x)g(v)dv.
\end{eqnarray*}
Thus,
$$\sup_{t_1,t_2}H(t_1,t_2)\leq \Gamma_x(\lambda_1,\lambda_2).$$
Studying now the function $H$, it is easily seen that
$$\infty=\sup_{t_1,t_2}H(t_1,t_2)\leq \Gamma_x(\lambda_1,\lambda_2)$$
whenever $\lambda_2/\lambda_1>v_1(x)$. Similarly, whenever
$\lambda_2/\lambda_1<v_0(x)$ , we obtain
$\infty=\Gamma_x(\lambda_1,\lambda_2)$.\hfill$\Box$
\vspace{5mm}

\noindent{\bf Proof of Corollary \ref{cor2}}\ \ The proof follows
straightforwardly from Theorem \ref{thm1} by making use of the
contraction principle with the following continuous function
\begin{eqnarray}\label{H}
H&:&\mathbb R_+^*\times\mathbb R\to\mathbb R\nonumber\\
&&(\lambda_1,\lambda_2)\to \frac{\lambda_2}{\lambda_1}.
\end{eqnarray}
Consequently,  $\hat r_n^l$ satisfies the LDP with the speed
$n\phi(h)$ and the rate function
\begin{eqnarray*}
  \gamma_x(\lambda)&:=&\inf \left\{\Gamma_x(\lambda_1,\lambda_2) : H(\lambda_1,\lambda_2)=\lambda\right\}
   = \inf \left\{\Gamma_x(\lambda_1,\lambda_2) : \lambda_2/\lambda_1=\lambda\right\} \\
   &=& \inf \left\{\Gamma_x(\lambda_1,\lambda_1\times\lambda)
   :\lambda_1>0\right\}.
\end{eqnarray*}
\hfill$\Box$\vspace{3mm}

\noindent{\bf Proof of Theorem \ref{thm2}}\ \ First of all, since
the rate function $\gamma_x$ is continuous, it follows by the
contraction principle, used with the continuous map $y\rightarrow
y-r^l(x)$, that for any $\lambda>0$,
\begin{equation}\label{cher1}
    \lim_{n\rightarrow\infty}\frac{1}{n\phi(h)}\log\mathbb{P}(|\hat{r}_n^l(x)-r^l(x)|>\lambda)=-\beta(x,\lambda).
\end{equation}
To state the uniform lower bound, it suffices to notice that for
any $x\in {\cal C}$, we have
\begin{eqnarray*}
\liminf_{n\rightarrow\infty}\frac{1}{n\phi(h)}\log\mathbb{P}(\|\hat{r}_n^l-r^l\|>\lambda)&\geq&
\liminf_{n\rightarrow\infty}\frac{1}{n\phi(h)}\log\mathbb{P}(|\hat{r}_n^l(x)-r^l(x)|>\lambda)\\
&\geq&-\beta(x,\lambda).
\end{eqnarray*}
Therefore,
\begin{equation}\label{lbound}
\liminf_{n\rightarrow\infty}\frac{1}{n\phi(h)}\log\mathbb{P}(\|\hat{r}_n^l-r^l\|>\lambda)\geq
-\rho(\lambda).
\end{equation}

\noindent Towards establishing the upper bound, observe first that
\begin{eqnarray}\label{decomp1}
  \|\hat{r}_n^l-r^l\|_{\cal C} &=& \max_{1\leq j\leq {\cal N}(\nu,{\cal C},d)}\sup_{x\in B_d(c_j,\nu)}|\hat{r}_n^l(x)-r^l(x)|,\nonumber\\
&\leq& \max_{1\leq j\leq {\cal N}(\nu,{\cal C},d)}|\hat{r}_n^l(c_j)-r^l(c_j)|+
\max_{1\leq j\leq {\cal N}(\nu,{\cal C},d)}\sup_{x\in B_d(c_j,\nu)}|\hat{r}_n^l(c_j)-\hat{r}_n^l(x)|\nonumber\\
&&+ \max_{1\leq j\leq {\cal N}(\nu,{\cal C},d)}\sup_{x\in
B_d(c_j,\nu)}|r^l(c_j)-r^l(x)|,
\end{eqnarray}
where $B_d(c_j,\nu)=\{x\in {\cal E} : d(c_j,x)\leq \nu\}$.
Assuming that $r^l$ is uniformly continuous on ${\cal C}$, it
follows, for any $\varepsilon>0$, that there exists $\nu>0$ such
that
$$\max_{1\leq j\leq {\cal N}(\nu,{\cal C},d)}\sup_{x\in
B_d(c_j,\nu)}|r^l(c_j)-r^l(x)|<\varepsilon.$$
For any $x\in {\cal C}$ and any $c\in {\cal C}$, observe that
\begin{eqnarray*}
 \hat{r}_n^l(x)- \hat{r}_n^l(c) &=& \frac{\sum_{i=1}^nl(Y_i)\left(K\left(\frac{d(x,X_i)}{h}\right)-K\left(\frac{d(c,X_i)}{h}\right)\right)}
 {\sum_{i=1}^nK\left(\frac{d(x,X_i)}{h}\right)} \\
  && +\sum_{i=1}^nl(Y_i)K\left(\frac{d(x,X_i)}{h}\right)\left[\frac{\sum_{i=1}^n\left(K\left(\frac{d(x,X_i)}{h}\right)-K\left(\frac{d(c,X_i)}{h}\right)\right)}
  {\sum_{i=1}^nK\left(\frac{d(x,X_i)}{h}\right)\sum_{i=1}^nK\left(\frac{d(c,X_i)}{h}\right)}\right].
\end{eqnarray*}
Assuming the kernel $K$ to be a Lipschitz function, it follows
that
$$\left|K\left(\frac{d(x,X_i)}{h}\right)-K\left(\frac{d(c,X_i)}{h}\right)\right|\leq
\frac{M}{h}|d(x,X_i)-d(c,X_i)|\leq \frac{M}{h}d(x,c),$$ where $M$
is a positive constant. Therefore, whenever the kernel $K$ is
bounded away from below by $K_0>0$, we have
\begin{eqnarray*}
  \sup_{x\in B_d(c,\nu)}|\hat{r}_n^l(x)- \hat{r}_n^l(c)| &\leq& \frac{M\nu}{K_0nh}\sum_{i=1}^n|l(Y_i)|+
  \frac{M\nu}{K_0^2nh}\sum_{i=1}^n|l(Y_i)|K\left(\frac{d(c,X_i)}{h}\right)\\
  &\leq&\frac{2M\nu}{K_0^2nh}\sum_{i=1}^n|l(Y_i)|K\left(\frac{d(c,X_i)}{h}\right).
\end{eqnarray*}
Thus, for any $\varepsilon>0$, by Markov's inequality, we obtain
\begin{eqnarray*}
\mathbb{P}\left(\sup_{x\in B_d(c,\nu)}|\hat{r}_n^l(x)-
\hat{r}_n^l(c)|>\varepsilon\right)&\leq&
\mathbb{P}\left(\sum_{i=1}^n|l(Y_i)|K\left(\frac{d(c,X_i)}{h}\right)>\varepsilon\frac{K_0^2nh}{2M\nu}\right)\\
&\leq&\frac{2M\nu}{\varepsilon
K_0^2h}\mathbb{E}\left(|l(Y_1)|K\left(\frac{d(c,X_1)}{h}\right)\right)
\end{eqnarray*}
Proceeding now similarly as in the proof of Theorem \ref{thm1}, we
obtain
\begin{eqnarray*}
\mathbb{E}\left(|l(Y_1)|K\left(\frac{d(c,X_1)}{h}\right)\right)\!\!\!&=&\phi(h)\left(\int\left[|l(v)|K(1)(f_v(c)+o(1))\right.\right.\\
&-&\!\!\!\!\left.\int_0^1\!\!\left.|l(v)|K'(u)(f_v(c)\!+\!o(1))(\tau_0(u)\!+\!o(1))du\right]g(v)dv\right)\\
&:=& \phi(h)A,
\end{eqnarray*}
where, by conditions (A4), A is a finite constant. Consequently,
there exists a positive constant $C$ such that
\begin{eqnarray*}
\mathbb{P}\left(\max_{1\leq j\leq {\cal N}(\nu,{\cal
C},d)}\sup_{x\in B_d(c_j,\nu)}|\hat{r}_n^l(x)-
\hat{r}_n^l(c_j)|>\varepsilon\right)&\leq&C{\cal N}(\nu,{\cal
C},d)\phi(h)\frac{\nu}{h\varepsilon}.
\end{eqnarray*}
Considering the decomposition in  the statement (\ref{decomp1}),
it is easily seen that
\begin{eqnarray*}
\mathbb{P}\left( \|\hat{r}_n^l-r^l\|_{\cal C}>\lambda\right)
&=&\mathbb{P}\left(\max_{1\leq j\leq {\cal N}(\nu,{\cal
C},d)}|\hat{r}_n^l(c_j)-r^l(c_j)|>\lambda-2\varepsilon\right)\\
&&\!\!\!\!\!\!\times\left(1+\frac{\mathbb{P}\left(\displaystyle{\max_{1\leq
j\leq {\cal N}(\nu,{\cal C},d)}\sup_{x\in
B_d(c_j,\nu)}}|\hat{r}_n^l(x)-
\hat{r}_n^l(c_j)|>\varepsilon\right)}{\mathbb{P}\left(\displaystyle{
\max_{1\leq j\leq {\cal N}(\nu,{\cal
C},d)}}|\hat{r}_n^l(c_j)-r^l(c_j)|>\lambda-2\varepsilon\right)}\right)\\
&\leq&\mathbb{P}\left( \max_{1\leq j\leq {\cal N}(\nu,{\cal
C},d)}|\hat{r}_n^l(c_j)-r^l(c_j)|>\lambda-2\varepsilon\right)\\
&&\times\left(1+\frac{C{\cal N}(\nu,{\cal
C},d)\phi(h)\frac{\nu}{h\varepsilon}}{\mathbb{P}\left(
\displaystyle{\max_{1\leq j\leq {\cal N}(\nu,{\cal
C},d)}}|\hat{r}_n^l(c_j)-r^l(c_j)|>\lambda-2\varepsilon\right)}\right)
\end{eqnarray*}
Since, for any $u>0$, $\log(1+u)\leq u$, it is obvious then that
\begin{eqnarray*}
\frac{1}{n\phi(h)}\log\mathbb{P}\left( \|\hat{r}_n^l-r^l\|_{\cal
C}\!>\!\lambda\right)\!\!&\!\!\leq\!\!&\!\!\frac{1}{n\phi(h)}\log\mathbb{P}\left(
\displaystyle{\max_{1\leq j\leq {\cal N}(\nu,{\cal
C},d)}}|\hat{r}_n^l(c_j)-r^l(c_j)|\!>\!\lambda-2\varepsilon\!\right)\\
&&+\frac{C{\cal N}(\nu,{\cal
C},d)\nu}{nh\varepsilon\mathbb{P}\left( \displaystyle{\max_{1\leq
j\leq {\cal N}(\nu,{\cal
C},d)}}|\hat{r}_n^l(c_j)-r^l(c_j)|>\lambda-2\varepsilon\right)}
\end{eqnarray*}
Therefore, from the statement (\ref{lbound}), we obtain
\begin{eqnarray*}
\lefteqn{\hspace*{-1.5cm}\frac{1}{n\phi(h)}\log\mathbb{P}\left(
\|\hat{r}_n^l-r^l\|_{\cal C}>\lambda\right)}\\
&\leq&\frac{1}{n\phi(h)}\log\mathbb{P}\left( \max_{1\leq j\leq
{\cal N}(\nu,{\cal
C},d)}|\hat{r}_n^l(c_j)-r^l(c_j)|>\lambda-2\varepsilon\right)\\
&&+\ \frac{C{\cal N}(\nu,{\cal
C},d)\nu}{nh\varepsilon}\exp\left\{n\phi(h)(\rho(\lambda)+o(1))\right\}\\
&\leq&\frac{1}{n\phi(h)}\log{\cal N}(\nu,{\cal C},d)\\
&&+\ \sup_{x\in{\cal C}}\frac{1}{n\phi(h)}\log\mathbb{P}\left(
\max_{1\leq j\leq {\cal N}(\nu,{\cal
C},d)}|\hat{r}_n^l(x)-r^l(x)|>\lambda-2\varepsilon\right)\\
&&+\ \frac{C{\cal N}(\nu,{\cal
C},d)\nu}{nh\varepsilon}\exp\left\{n\phi(h)(\rho(\lambda)+o(1))\right\}.
\end{eqnarray*}
Making use of the condition (\ref{entropy}), it follows that
$$\limsup_{n\rightarrow\infty}\frac{1}{n\phi(h)}\log\mathbb{P}\left( \|\hat{r}_n^l-r^l\|_{\cal
C}>\lambda\right)\leq -\rho(\lambda-2\varepsilon).$$ The proof is
achieved while making $\varepsilon$ tend to zero since the
function $\rho(\lambda)$ is continuous. \hfill$\Box$ \vspace{5mm}

\noindent{\bf Proof of Proposition \ref{prop1}}\ \ Observe first,
for any positive real numbers $\lambda$ and $\lambda_1$, that
\begin{eqnarray*}
  \rho(\lambda) &=& \inf_{x\in\C}\beta(x,\lambda)= \inf_{x\in\C}\left\{\beta(x,\lambda)-\beta(x,\lambda_1)+\beta(x,\lambda_1)\right\}\\
 &\leq& \inf_{x\in\C}\beta(x,\lambda_1)+
 \sup_{x\in\C}\{\beta(x,\lambda)-\beta(x,\lambda_1)\}.
\end{eqnarray*}
Therefore, we have
$$|\rho(\lambda)-\rho(\lambda_1)|\leq
\sup_{x\in\C}|\beta(x,\lambda)-\beta(x,\lambda_1)|.$$ Considering
the statement (\ref{beta}), we obtain for $\lambda<\lambda_1$
\begin{eqnarray*}
|\rho(\lambda)-\rho(\lambda_1)|&\leq&
\sup_{x\in\C}\left|\inf_{|\alpha|\geq\lambda}\gamma_x(\alpha+r^l(x))-\inf_{|\alpha|\geq\lambda_1}\gamma_x(\alpha+r^l(x))\right|\\
&\leq&\sup_{x\in\C}\left|\inf_{|\alpha|\geq\lambda}(\gamma_x(\alpha+r^l(x))-\gamma_x(\lambda+r^l(x)))\right.\\
&&\ \ \ \ \left.-\inf_{|\alpha|\geq\lambda_1}(\gamma_x(\alpha+r^l(x))-\gamma_x(\lambda+r^l(x)))\right|\\
&\leq&\!\!\!\sup_{x\in\C}\left\{\inf_{|\alpha|\geq\lambda}|\gamma_x(\alpha+r^l(x))-\gamma_x(\lambda+r^l(x)))|\right.\\
&&\ \ \ \ \left.+
\inf_{|\alpha|\geq\lambda_1}|\gamma_x(\alpha+r^l(x))-\gamma_x(\lambda+r^l(x)))|\right\}\\
&\leq&\!\!\!2\sup_{x\in\C}\left|\gamma_x(\lambda+r^l(x))-\gamma_x(\lambda_1+r^l(x)))\right|.
\end{eqnarray*}
Taking into account the shape of the rate function $\gamma_x$
given in Corollary \ref{cor2}, we obtain
\begin{eqnarray*}
|\rho(\lambda)-\rho(\lambda_1)|&\leq&
2\sup_{x\in\C}\left|\inf_{\delta}\Gamma_x(\delta,\delta(\lambda+r^l(x)))-\inf_{\delta}\Gamma_x(\delta,\delta(\lambda_1+r^l(x)))\right|\\
&\leq&
2\sup_{x\in\C}\left|\inf_{\delta}(\Gamma_x(\delta,\delta(\lambda+r^l(x)))-\Gamma_x(\lambda,\lambda(\lambda_1+r^l(x))))\right.\\
&&\ \ \ \ \ -\left.\inf_{\delta}(\Gamma_x(\delta,\delta(\lambda_1+r^l(x)))-\Gamma_x(\lambda,\lambda(\lambda_1+r^l(x))))\right|\\
&\leq&
2\sup_{x\in\C}\left|\Gamma_x(\lambda,\lambda(\lambda+r^l(x)))-\Gamma_x(\lambda,\lambda(\lambda_1+r^l(x)))\right|\\
&&+2\sup_{x\in\C}\left|\Gamma_x(\lambda_1,\lambda_1(\lambda_1+r^l(x)))-\Gamma_x(\lambda,\lambda(\lambda_1+r^l(x)))\right|.
\end{eqnarray*}
It suffices now to use the fact that
$\{\Gamma_x(\lambda,\mu)\}_{x\in\C}$ is a equi-continuous family
of functions and that the regression function $r^l$ is bounded to
achieve the proof. \hfill$\Box$ \vspace{5mm}

\noindent{\bf Proof of Proposition \ref{prop2}}\ \ \noindent
Observe by the condition (A3)(i)-(iii) that $\Phi^x$ is a
differentiable function and that
$$\Gamma_x(\lambda,\mu)=\lambda s_x+\mu t_x-\Phi^x(s_x,t_x),$$
where $s_x$ and $t_x$ are solutions of equations
\begin{eqnarray}\label{partial}
\mbox{(i)}\ \ \frac{\partial \Phi^x}{\partial s}(s,t)=\lambda,\ \
\ \mbox{(ii)}\ \ \frac{\partial \Phi^x}{\partial t}(s,t)=\mu.
\end{eqnarray}
respectively. Therefore, for any $\lambda_1$ and $\mu_1$, we have
\begin{eqnarray*}
\Gamma_x(\lambda,\mu)&=&\lambda_1 s_x+\mu_1 t_x-\Phi^x(s_x,t_x)
+(\lambda-\lambda_1)s_x+(\mu-\mu_1)t_x\\
&\leq& \sup_{s,t}\{\lambda_1 s +\mu_1 t-\Phi^x(s,t)\}
+|(\lambda-\lambda_1)s_x|+|(\mu-\mu_1)t_x|\\
&\leq&\Gamma_x(\lambda_1,\mu_1)+|(\lambda-\lambda_1)s_x|+|(\mu-\mu_1)t_x|
\end{eqnarray*}
and then
$$|\Gamma_x(\lambda,\mu)-\Gamma_x(\lambda_1,\mu_1)|\leq|\lambda-\lambda_1||s_x|+|\mu-\mu_1||t_x|.$$
We have now to state that both $\displaystyle{\sup_{x\in\C}|t_x|}$
and $\displaystyle{\sup_{x\in\C}|t_x|}$ are finite.

Observe from the statement (\ref{partial})(i), whenever $t_x\geq
0$,  that
\begin{eqnarray*}
  \lambda &=& e^{s_x}\left(\int_{A_l} e^{t_x+t_x(l(v)-1)}f_v(x)g(v)dv+ \int_{\bar{A}_l} e^{t_x+t_x(l(v)-1)}f_v(x)g(v)dv\right)\\
  &\geq&e^{s_x+t_x}\int_{\bar{A}_l} e^{t_x(l(v)-1)}f_v(x)g(v)dv \geq e^{s_x+t_x}\int_{\bar{A}_l} d(v)g(v)dv.
\end{eqnarray*}
Therefore, we have
$$ s_x+t_x\leq \log\left(\frac{\lambda}{\int_{\bar{A}_l}
 d(v)g(v)dv}\right).$$
%
Whenever $t_x<0$, similarly, we obtain
$$ s_x+t_x\leq \log\left(\frac{\lambda}{\int_{A_l}
 d(v)g(v)dv}\right).$$
Therefore, whatever the value that may take $t_x$, we obtain
\begin{equation}\label{S+T3}
 s_x+t_x\leq \max\left\{\log\left(\frac{\lambda}{\int_{\bar{A}_l}
 d(v)g(v)dv}\right),\log\left(\frac{\lambda}{\int_{A_l}
 d(v)g(v)dv}\right)\right\}.
\end{equation}
On another hand, whenever $t_x\geq 0$, we have
\begin{eqnarray*}
  \lambda &=& e^{s_x}\left(\int_{B_l} e^{-t_x+t_x(l(v)+1)}f_v(x)g(v)dv+ \int_{\bar{B}_l} e^{-t_x+t_x(l(v)+1)}f_v(x)g(v)dv\right)\\
  &\geq&e^{s_x-t_x}\int_{\bar{B}_l} e^{t_x(l(v)+1)}f_v(x)g(v)dv \geq e^{s_x-t_x}\int_{\bar{B}_l} d(v)g(v)dv.
\end{eqnarray*}
Consequently, we have
$$ s_x-t_x\leq \log\left(\frac{\lambda}{\int_{\bar{B}_l}
 d(v)g(v)dv}\right).$$
Whenever $t_x<0$, similarly, we obtain
 $$s_x-t_x\leq \log\left(\frac{\lambda}{\int_{B_l}
 d(v)g(v)dv}\right).$$
Therefore, for any value that may take $t_x$, we obtain
\begin{equation}\label{S+T6}
 s_x-t_x\leq \max\left\{\log\left(\frac{\lambda}{\int_{\bar{B}_l}
 d(v)g(v)dv}\right),\log\left(\frac{\lambda}{\int_{B_l}
 d(v)g(v)dv}\right)\right\}.
\end{equation}
Considering the statements (\ref{S+T3}) and (\ref{S+T6}), it
follows that
\begin{eqnarray}\label{s1}
  s_x &\leq& \frac{1}{2}\max\left\{\log\left(\frac{\lambda}{\int_{\bar{A}_l}
 d(v)g(v)dv}\right),\log\left(\frac{\lambda}{\int_{A_l}
 d(v)g(v)dv}\right)\right\}\nonumber\\
  && \ \ \ +\ \ \frac{1}{2}\max\left\{\log\left(\frac{\lambda}{\int_{\bar{B}_l}
 d(v)g(v)dv}\right),\log\left(\frac{\lambda}{\int_{B_l}
 d(v)g(v)dv}\right)\right\}=:s_1.
\end{eqnarray}

Taking into consideration the definition of the function $V_x$
given in the statement (\ref{defV}), it is clear, for any
$t\in\R$, that
$$V_d(t):=\frac{\int
e^{tl(v)}l(v)(D(v)\mathds{1}_{C_l}(v)+d(v)\mathds{1}_{\bar{C}_l}(v))g(v)dv}{\int
e^{tl(v)}D(v)g(v)dv}\leq V_x(t).$$ Consequently, considering the
inverses of the functions $V_d$ and $V_x$ together with equations
of the statement (\ref{partial}), since $V_x$ is a nondecreasing
function, it follows, whenever
$\displaystyle{\frac{\mu}{\lambda}>v_{0,d}:=\inf_tV_d(t)}$, that
\begin{equation}\label{bsup}
t_1:=V_d^{-1}\left(\frac{\mu}{\lambda}\right):=\inf\left\{s :
V_d(s)\geq \frac{\mu}{\lambda}\right\} \geq \inf\left\{s :
V_x(s)\geq \frac{\mu}{\lambda}\right\}=t_x.
\end{equation}
Observe now, for any $t\in\R$, that we have
$$V_D(t):=\frac{\int
e^{tl(v)}l(v)(d(v)\mathds{1}_{C_l}(v)+D(v)\mathds{1}_{\bar{C}_l}(v))g(v)dv}{\int
e^{tl(v)}d(v)g(v)dv}\geq V_x(t).$$
Similarly as above, whenever
$\displaystyle{\frac{\mu}{\lambda}<v_{1,D}:=\sup_tV_D(t)}$, we
have
\begin{equation}\label{binf}
t_0:=V_D^{-1}\left(\frac{\mu}{\lambda}\right):=\inf\left\{s :
V_D(s)\geq \frac{\mu}{\lambda}\right\} \leq \inf\left\{s :
V_x(s)\geq \frac{\mu}{\lambda}\right\}=t_x.
\end{equation}

Moreover, it is obvious from the statement (\ref{partial})(i) that
$$
 e^{s_x}\left(\int_{C_l} e^{t_0l(v)}D(v)g(v)dv+\int_{\bar{C}_l}
 e^{t_1l(v)}D(v))g(v)dv\right)\geq \lambda.
$$
Therefore, for any $\lambda>0$, we have
\begin{equation}\label{s0}
s_x\geq\log\left(\frac{\lambda}{\int_{C_l}
e^{t_0l(v)}D(v)g(v)dv+\int_{\bar{C}_l}
 e^{t_1l(v)}D(v))g(v)dv}\right)=:s_0.
\end{equation}
Consequently, making use of the statements (\ref{s1}),
(\ref{bsup}), (\ref{binf})  and (\ref{s0}), it is clear that there
exist  finite numbers $S_0$ and $T_0$ such that
$$\sup_{x\in\C}|\Gamma_x(\lambda,\mu)-\Gamma_x(\lambda_1,\mu_1)|\leq|\lambda-\lambda_1||S_0|+|\mu-\mu_1||T_0|.$$
This establishes that the family of functions
$\{\Gamma_x(\lambda,\mu)\}_{x\in\C}$ is a equi-continuous and
achieves the proof.$\hfill\Box$\vspace{5mm}

\centerline{\large REFERENCES}

\begin{enumerate}

\item\textsc{Bosq, D.} (1998). { Nonparametric Statistics for Stochastic
Processes}. Lecture Note in Statistics. Springer, New York.

\item\textsc{Dembo,  A. \& Zeitouni, O.} (1998). \textit{Large deviations
techniques and applications}. Second edition, Springer-Verlag, New
york.

\item\textsc{Deuschel, J.D. \& Stroock, D.W.} (1989). \textit{Large
deviations.} Academic Press.

\item\textsc{Ezzahrioui, M. and Ould-Sa\"{i}d, E.} (2008). Asymptotic
normality of a nonparametric estimator of the conditional mode
function for functional data. {\it J. Nonparametric. Statist.},
{\bf 20},  3--18.

\item\textsc{Ferraty, F. and Vieu, P.} (2000). Dimension fractale et
estimation de la r\'egression dans des espaces vectoriels
semi-norm\'es. {\it C. R. Acad. Sci. Paris S\'er. I Math.},  {\bf
330},  139--142.

\item\textsc{Ferraty, F. and Vieu, P.} (2004). Nonparametric models for
functional data, with applications in regression, time series
prediction and curve discrimination. The International Conference
on Recent Trends and Directions in Nonparametric Statistics. {\it
J. Nonparametric Statist.},  {\bf 16}, 111--125.

\item\textsc{Ferraty, F., Laksaci, A. and Vieu, P.} (2006).  Estimating some
characteristics of the conditional distribution in nonparametric
functional models. {\it Stat. Inference Stoch. Process.}, {\bf 9},
47-76.

\item\textsc{Ferraty, F. and Vieu, P.} (2006). {\it Nonparametric
functional data analysis. Theory and practice.} Springer Series in
Statistics. Springer, New York

\item\textsc{Ferraty, F., Mas, A. and Vieu, P.} (2007). Nonparametric
regression of functional data: inference and practical aspects.
{\it Aust.N.Z.J.Statist.}, {\bf 49}, 267-286.

\item\textsc{Gasser, T., Hall, P. and Presnell, B.}  (1998). Nonparametric
estimation of the mode of a distribution of random curves. {\it J.
Roy. Statist. Soc.  Ser B}, {\bf 60}, 681-691.

\item\textsc{La\"ib, N. and Louani, D.} (2010). Nonparametric Kernel Regression Estimation for  Functional Stationary Ergodic  Data:
 Asymptotic Properties. \textit{J. Multivariate Analysis}, {\bf
 101}, 2266-2281.

\item\textsc{Louani, D.} (1999). Some large deviations limit theorems in conditional nonparametric statistics.
\textit{Statistics}, \textbf{33}, 171-196.

\item\textsc{Louani, D. and Ould Maouloud, S. M.} (2011). Some Functional Large Deviations Principles in Nonparametric Function
Estimation. \textit{J. Theoretical Probability}, \textbf{24}. In
press.

\item\textsc{Nikitin, Ya.} (1995). \textit{Asymptotic efficiency of
non-parametric tests.} Cambridge University Press, Cambridge.

\item\textsc{Ould Maouloud, S. M.} (2008). Some uniform large deviation results
in nonparametric function estimation. \textit{J. Nonparametric
Statist.}, \textbf{20}, 129 -152.

\item\textsc{van der Vaart, A. W. \& Wellner, J. A.} (1996). Weak convergence and
empirical processes. With applications to statistics. Springer
Series in Statistics. Springer-Verlag, New York.

\end{enumerate}

\end{document}